\documentclass[a4paper,11pt,twoside,reqno]{amsart}

\usepackage[utf8]{inputenc}
\usepackage[plainpages=false,pdfpagelabels=true]{hyperref}
\usepackage{amssymb,amsthm,esint}
\usepackage[margin=1in]{geometry}

\numberwithin{equation}{section}

\newtheorem{Satz}{Theorem}[section]
\newtheorem{Prop}[Satz]{Proposition}
\newtheorem{Lem}[Satz]{Lemma}

\theoremstyle{definition}
\newtheorem{Dfn}[Satz]{Definition}
\newtheorem{Bem}[Satz]{Remark}

\newcommand{\vol}{\operatorname{vol}}
\newcommand{\cL}{{\mathcal{L}}}
\newcommand{\tr}{\operatorname{Tr}}
\newcommand{\sff}{\mathrm{I\!I}}
\newcommand{\N}{\ensuremath{\mathbb{N}}}
\newcommand{\R}{\ensuremath{\mathbb{R}}}

\theoremstyle{remark}

\newcommand{\dv}{\text{ }dV}
\allowdisplaybreaks[1]

\renewcommand{\epsilon}{\varepsilon}

\title[Some analytic results on interpolating sesqui-harmonic maps]{Some analytic results on interpolating sesqui-harmonic maps}
\author{Volker Branding}
\date{\today}
\address{University of Vienna, Faculty of Mathematics\\
Oskar-Morgenstern-Platz 1, 1090 Vienna, Austria\\}
\email{volker.branding@univie.ac.at}

\begin{document}

\begin{abstract}
In this article we study various analytic aspects of interpolating sesqui-harmonic maps between Riemannian manifolds where 
we mostly focus on the case of a spherical target.
The latter are critical points of an energy functional that interpolates between the functionals for harmonic and biharmonic maps.
In the case of a spherical target we will derive a conservation law and use it to show the smoothness of weak solutions.
Moreover, we will obtain several classification results for interpolating sesqui-harmonic maps.
\end{abstract} 
 
\maketitle

\section{Introduction and results}
Harmonic maps are among the most important variational problems in geometry, analysis and physics.
Given a map \(\phi\colon M\to N\) between two Riemannian manifolds \((M,g)\) and \((N,h)\)
they are defined as critical points of the Dirichlet energy
\begin{align}
\label{harmonic-energy}
E(\phi)=\int_M|d\phi|^2\dv.
\end{align}
The first variation of \eqref{harmonic-energy} is characterized by the vanishing of the so-called
\emph{tension field} which is given by
\begin{align}
\label{harmonic-map-equation}
0=\tau(\phi):=\tr_g\bar\nabla d\phi.
\end{align}
Here, \(\bar\nabla\) represents the connection on \(\phi^\ast TN\).
Solutions of \eqref{harmonic-map-equation} are called \emph{harmonic maps}.
The harmonic map equation is a semilinear, elliptic second-order partial differential equation for which many results 
on existence and qualitative behavior of its solutions could be achieved over the years.

A higher order generalization of harmonic maps, that receives growing attention, are 
the so-called \emph{biharmonic maps}.
These arise as critical points of the bienergy for a map between two Riemannian manifolds which is given by
\begin{align*}
E_2(\phi)=\int_M|\tau(\phi)|^2\dv
\end{align*}
and are characterized by the vanishing of the \emph{bitension field}
\begin{align*}
0=\tau_2(\phi):=\bar\Delta\tau(\phi)-R^N(d\phi(e_j),\tau(\phi))d\phi(e_j).
\end{align*}
Here, \(\bar\Delta\) is the connection Laplacian on \(\phi^\ast TN\), \(\{e_j\},j=1,\ldots,m=\dim M\) an orthonormal basis of \(TM\) and 
\(R^N\) denotes the curvature tensor of the target manifold \(N\). Moreover, we apply the Einstein summation convention, meaning
that we sum over repeated indices. 

In contrast to the harmonic map equation the biharmonic map equation is a semilinear elliptic equation of fourth order 
such that its study comes with additional difficulties.

It can be directly seen that every harmonic map is also biharmonic. However,
a biharmonic map can be non-harmonic in which case it is called \emph{proper biharmonic}.
Many conditions, both of analytic and geometric nature, are known that force a biharmonic map to be harmonic,
see for example \cite{MR3834926,MR4040175} and references therein for a recent overview.

An extensive study of higher order energy functionals for maps between Riemannian manifolds
was recently initiated in \cite{branding2019higher}.

In this article we want to focus on the critical points of an energy functional
that interpolates between the energy functionals for harmonic and biharmonic maps which is given by
\begin{align}
\label{energy-functional}
E_{\delta_1,\delta_2}(\phi)=\delta_1\int_M|d\phi|^2\dv+\delta_2\int_M|\tau(\phi)|^2\dv
\end{align}
with \(\delta_1, \delta_2\in\R\). Various versions of this functional had already been studied by a number of
mathematicians, a general study of \eqref{energy-functional} was recently initiated by the author \cite{MR4058514}.

The critical points of \eqref{energy-functional} will be referred to as \emph{interpolating sesqui-harmonic maps} and are given by
\begin{align}
\label{euler-lagrange}
\delta_2\bar\Delta\tau(\phi)=\delta_2 R^N(d\phi(e_j),\tau(\phi))d\phi(e_j)+\delta_1\tau(\phi).
\end{align}
For more background on both harmonic and biharmonic maps and references that study versions of \eqref{energy-functional}
we refer to the introduction of \cite{MR4058514}.

So far, only few results on interpolating sesqui-harmonic maps have been established.
Besides a number of general features \cite{MR4058514},
the unique continuation property for interpolating sesqui-harmonic maps has been proved in \cite{MR3990379}.

In this article we will mostly focus on analytic aspects of interpolating sesqui-harmonic maps
building on the regularity theory developed for biharmonic maps \cite{MR1692148,MR2375314}.
The regularity for a closely related problem has been studied in \cite{MR1809291,MR2094552}.

We will establish the regularity of weak solutions in the case of a spherical target.
In order to achieve this result we will exploit the large amount of symmetry of the sphere
building on some earlier work for harmonic maps to target manifolds with a sufficient amount of symmetry \cite{MR1085633}
and biharmonic maps to spheres \cite{MR1692148}.
We will not deduce the regularity from the equation for interpolating sesqui-harmonic maps itself.
At first, we will derive a conserved quantity, arising due to the symmetries of the sphere.
In a second step, using the conserved quantity, we will establish the regularity of weak solutions by applying
the regularity theory for biharmonic maps to spheres \cite{MR2094320}.
We hope that the framework used in this article may turn out to be useful for other geometric 
variational problems of higher order in the future.

After the study of the regularity of weak solutions for a spherical target
we will also provide some remarks on interpolating sesqui-harmonic immersions to spheres.

Finally, we will prove a classification result for solutions of \eqref{euler-lagrange}
in the case of a Euclidean domain. This result shows that under assumption on energy, dimension
and the signs of both \(\delta_1,\delta_2\) solutions of \eqref{euler-lagrange} have
to be harmonic maps or even trivial.

Aiming in a similar direction we will also discuss how to obtain a monotonicity type formula
for solutions of \eqref{euler-lagrange} where we again, for simplicity, focus on the case of a Euclidean domain.

Throughout this paper we will make use of the following conventions:
Whenever we will make use of indices, we will use
Latin indices \(i,j,k\) for indices on the domain ranging from \(1\) to \(m=\dim M\)
and Greek indices \(\alpha,\beta,\gamma\) for indices on the target
which take values between \(1\) and \(n=\dim N\). In addition, local coordinates on the domain
will be denoted by \(x^i\) and for local coordinates on the target
we will use \(y^\alpha\).

In this article the curvature tensor is defined as \(R(X,Y)Z:=[\nabla_X,\nabla_Y]Z-\nabla_{[X,Y]}Z\)
such that the sectional curvature is given by \(K(X,Y)=R(X,Y,Y,X)\).
For the Laplacian acting on functions \(f\in C^\infty(M)\) we choose the convention
\(\Delta f=\operatorname{div}\operatorname{grad}f\),
for sections in the vector bundle \(\phi^{\ast}TN\) we make
the choice \(\bar\Delta=\tr(\nabla^{\phi^\ast TN}\nabla^{\phi^\ast TN})\).

This article is organized as follows:
In Section 2 we study the regularity of weak interpolating sesqui-harmonic maps 
for a spherical target and make some comments on interpolating sesqui-harmonic immersions to spheres.
In Section 3 we first prove a classification result for interpolating sesqui-harmonic maps 
from Euclidean space and also give a monotonicity type formula.

\section{Interpolating sesqui-harmonic maps to spheres}
In this section we study several aspects of \eqref{euler-lagrange} in the case of a spherical target.
We make use of the inclusion map \(\iota\colon S^n\to\R^{n+1}\) and consider the composite map
\(\varphi:=\iota\circ\phi\colon M\to\R^{n+1}\).

If \(\varphi\colon M\to S^n\subset\R^{n+1}\) then \eqref{euler-lagrange} acquires the form 
\begin{align}
\label{euler-lagrange-sphere}
\delta_2(\Delta^2\varphi+(|\Delta\varphi|^2+\Delta|d\varphi|^2+2\langle d\varphi,\nabla\Delta\varphi\rangle+2|d\varphi|^4)\varphi+2\nabla(|d\varphi|^2d\varphi))=\delta_1(\Delta\varphi+|d\varphi|^2\varphi)
\end{align}
assuming that the sphere is equipped with the constant curvature metric.
For a derivation of \eqref{euler-lagrange-sphere} see \cite[Proposition 2.5]{MR4058514}.

First, let us make the following observation:
\begin{Bem}
Besides the energy functional for interpolating sesqui-harmonic maps, there exists a similar functional 
that has received a lot of attention. In the case that the target \(N\) is realized as a submanifold of some \(\R^q\)
such that \(\phi\colon M\to\R^q\) the functional for \emph{extrinsic biharmonic maps} is given by 
\begin{align*}
E_{ext}(\phi)=\int_M(|\tau(\phi)|^2+|\sff(d\phi,d\phi)|^2)\dv.
\end{align*}
Here, \(\sff\) denotes the second fundamental form of the embedding.

It is well known that in the case of a spherical target the tension field acquires the simple form
\begin{align*}
d\iota(\tau(\phi))=\Delta\varphi+|d\varphi|^2\varphi,
\end{align*}
where \(\phi\colon M\to S^n\), \(\iota\colon S^n\to\R^{n+1}\) and \(\varphi=\iota\circ\phi\colon M\to\R^{n+1}\).
Consequently, we obtain for \(\varphi\colon M\to S^n\subset\R^{n+1}\) that
\begin{align*}
E_{ext}(\varphi)=\int_M(|\tau(\varphi)|^2+|d\varphi|^4)\dv=\int_M|\Delta\varphi|^2\dv.
\end{align*}
Assuming that \(\delta_2>0\) we obtain the following inequality 
\begin{align*}
E_{\delta_1,\delta_2}(\varphi)
&\leq\int_M(\delta_2|\tau(\varphi)|^2+\delta_2|d\varphi|^4)\dv+\frac{\delta_1^2}{4\delta_2}\vol(M,g) \\
&\leq\delta_2 E_{ext}(\varphi)+\frac{\delta_1^2}{4\delta_2}\vol(M,g).
\end{align*}
Hence, under the assumption that \(\delta_2>0\) the energy functional for interpolating sesqui-harmonic maps 
is bounded from above by the energy functional for extrinsic biharmonic maps and a constant. 
Consequently, we should expect that the critical points of both functionals share common properties.
\end{Bem}

\subsection{Conserved currents}
In this section we will demonstrate how one can obtain a conservation law from \eqref{euler-lagrange-sphere}
by exploiting the symmetries of the sphere. As the energy functional for interpolating sesqui-harmonic maps \eqref{energy-functional}
is invariant under isometries of the target manifold we are getting a conserved quantity via Noether's theorem.
We will illustrate various methods how to explicitly compute this conserved quantity.

The next Lemma is similar to \cite[Lemma 2.2]{MR2094320}.

\begin{Lem}
Suppose that \(\varphi\colon M\to S^n\subset\R^{n+1}\) is a smooth solution of \eqref{euler-lagrange-sphere}.
Then the following conservation law holds true
\begin{align}
\label{conservation-law-sphere-a}
\nabla\big(\delta_2(\nabla\Delta\varphi\wedge\varphi-\Delta\varphi\wedge\nabla\varphi+2|d\varphi|^2\nabla\varphi\wedge\varphi)-\delta_1\nabla\varphi\wedge\varphi)\big)=0,
\end{align}
which can also be written in the following form
\begin{align}
\label{conservation-law-sphere-b}
\Delta (\nabla(\nabla\varphi\wedge\varphi))=\nabla\big(2\Delta\varphi\wedge\nabla\varphi-2|d\varphi|^2\nabla\varphi\wedge\varphi+\frac{\delta_1}{\delta_2}\nabla\varphi\wedge\varphi\big).
\end{align}
\end{Lem}

\begin{proof}
Wedging \eqref{euler-lagrange-sphere} with \(\varphi\) we get 
\begin{align*}
\delta_2\Delta^2\varphi\wedge\varphi+2\delta_2\nabla(|d\varphi|^2d\varphi)\wedge\varphi=\delta_1\Delta\varphi\wedge\varphi.
\end{align*}
By a direct calculation we find
\begin{align*}
\Delta^2\varphi\wedge\varphi=&\nabla(\nabla\Delta\varphi\wedge\varphi-\Delta\varphi\wedge\nabla\varphi), \\
\nabla(|d\varphi|^2\nabla\varphi)\wedge\varphi=&\nabla (|d\varphi|^2\nabla\varphi\wedge\varphi), \\
\Delta\varphi\wedge\varphi=&\nabla(\nabla\varphi\wedge\varphi),
\end{align*}
which proves the first assertion.
The second formula follows from the identity
\begin{align*}
\Delta^2\varphi\wedge\varphi=\Delta(\nabla(\nabla\varphi\wedge\varphi))-2\nabla(\Delta\varphi\wedge\nabla\varphi).
\end{align*}
\end{proof}

For harmonic maps to target spaces with a certain amount of symmetry such kinds of conservation laws have been obtained in
\cite{MR1078114}, see also \cite{MR3735550} for further applications.

As a next step we will show how the conservation law \eqref{conservation-law-sphere-b} 
can also be obtained in the case that we only have a weak solution of \eqref{euler-lagrange-sphere}.
A weak solution of \eqref{euler-lagrange-sphere} corresponds to \(\varphi\in W^{2,2}(M,S^n)\) which solves
\eqref{euler-lagrange-sphere} in a distributional sense.
To this end we recall the following facts:
\begin{Dfn}
A vector field \(X\) is called \emph{Killing vector field} on \((N,h)\) if
\[
\mathcal{L}_Xh=0,
\]
where \(\cL\) represents the Lie derivative of the metric. In terms of local coordinates we have
\begin{align}
\label{killing-definition}
0=(\mathcal{L}_Xh)_{\alpha\beta}=h_{\beta\gamma}\nabla_{\partial_{y^\alpha}}X^\gamma+h_{\alpha\gamma}\nabla_{\partial_{y^\beta}}X^\gamma.
\end{align}
\end{Dfn}

The group \(SO(n+1)\) acts isometrically on \(S^n\).
The set of Killing vector fields on \(S^n\) can be identified with the Lie algebra \(\mathfrak{so}(n+1)\) of \(SO(n+1)\).
In addition, \(\mathfrak{so}(n+1)\) can be represented as \((n+1)\times (n+1)\) skew-symmetric real-valued matrices.

\begin{Prop}
Suppose that \(\varphi\in W^{2,2}(M,S^n)\) is a weak solution of \eqref{euler-lagrange-sphere}
and let \(X\) be a Killing vector field on \(S^n\).
Then the following conservation law holds true
\begin{align}
\label{conservation-law-sphere-weak}
\int_M\big(\delta_2(\langle\Delta\varphi,\nabla X(\varphi)\rangle
-\langle\nabla\Delta\varphi,X(\varphi)\rangle
-2|d\varphi|^2\langle\nabla\varphi,X(\varphi)\rangle)+\delta_1\langle\nabla\varphi, X(\varphi)\rangle\big)\nabla\eta\dv
=0
\end{align}
for all \(\eta\in C^\infty(M)\).
\end{Prop}

\begin{proof}
We test \eqref{euler-lagrange-sphere} with \(X(\varphi)\eta\), where \(X(\varphi)\) is a Killing vector field on \(S^n\) and \(\eta\in C^\infty(M)\).
Then, we find
\begin{align*}
\delta_2\int_M\langle\Delta^2\varphi,X(\varphi)\rangle\eta\dv+2\delta_2\int_M\langle\nabla(|d\varphi|^2\nabla\varphi),X(\varphi)\rangle\eta\dv=\delta_1\int_M\langle\Delta\varphi, X(\varphi)\rangle\eta\dv,
\end{align*}
where we used that \(X(\varphi)\perp\varphi\).
Moreover, using integration by parts, we find
\begin{align*}
\int_M\langle\Delta\varphi, X(\varphi)\rangle\eta\dv=&-\int_M\underbrace{\langle\nabla\varphi,\nabla X(\varphi)\rangle}_{=0}\eta\dv-\int_M\langle\nabla\varphi, X(\varphi)\rangle\nabla\eta\dv,\\
\int_M\langle\nabla(|d\varphi|^2\nabla\varphi),X(\varphi)\rangle\eta\dv=&-\int_M|d\varphi|^2\underbrace{\langle\nabla\varphi,\nabla X(\varphi)\rangle}_{=0}\eta\dv
-\int_M|d\varphi|^2\langle\nabla\varphi,X(\varphi)\rangle\nabla\eta\dv.
\end{align*}
Here, we used that \(X(\varphi)\) is a solution of \eqref{killing-definition}.

Regarding the term that contains the Bi-Laplacian we find
\begin{align*}
\int_M\langle\Delta^2\varphi,X(\varphi)\rangle\eta\dv=&\int_M\langle\Delta\varphi,\Delta X(\varphi)\rangle\eta\dv+2\int_M\langle\Delta\varphi,\nabla X(\varphi)\rangle\nabla\eta\dv \\
&+\int_M\langle\Delta\varphi, X(\varphi)\rangle\Delta\eta\dv.
\end{align*}
To manipulate the first term on the right hand side we differentiate the equation for a Killing vector field on \(S^n\) \eqref{killing-definition} and obtain
\begin{align*}
\Delta X_\alpha(\varphi)+R_{\alpha\beta}X^\beta(\varphi)=0,
\end{align*}
where \(R_{\alpha\beta}\) denotes the Ricci curvature of \(S^n\).
Since the Ricci curvature on the sphere satisfies \(R_{\alpha\beta}=(n-1)h_{\alpha\beta}\) we find
\begin{align*}
\Delta X_\alpha(\varphi)=(1-n)X_\alpha(\varphi).
\end{align*}
Remember that \(0=\langle\varphi,X(\varphi)\rangle\), applying the Laplacian and using that \(X(\varphi)\) is a Killing vector field we obtain
\begin{align*}
0=\langle\Delta\varphi,X(\varphi)\rangle+2\underbrace{\langle\nabla\varphi,\nabla X(\varphi)\rangle}_{=0}+\langle\varphi,\Delta X(\varphi)\rangle.
\end{align*}
Hence, we find (in the sense of distributions)
\begin{align*}
\langle\Delta\varphi,\Delta X(\varphi)\rangle=(1-n)\langle\Delta\varphi,X(\varphi)\rangle=-(1-n)\langle\varphi,\Delta X(\varphi)\rangle=-(1-n)^2\langle\varphi,X(\varphi)\rangle=0.
\end{align*}
Consequently, we obtain
\begin{align*}
\int_M\langle\Delta^2\varphi,X(\varphi)\rangle\eta\dv=&2\int_M\langle\Delta\varphi,\nabla X(\varphi)\rangle\nabla\eta\dv
+\int_M\langle\Delta\varphi, X(\varphi)\rangle\Delta\eta\dv \\
=&\int_M(\langle\Delta\varphi,\nabla X(\varphi)\rangle-\langle\nabla\Delta\varphi, X(\varphi)\rangle)\nabla\eta\dv
\end{align*}
and the claim follows by combining the different equations.
\end{proof}

\begin{Bem}
Note that \eqref{conservation-law-sphere-weak} can be considered as the distributional version 
of \eqref{conservation-law-sphere-a}.
\end{Bem}

\begin{Bem}
There exists another way how to derive the conservation law \eqref{conservation-law-sphere-a}.
We recall the following result \cite[Proposition 2.11]{MR4058514}:

If \(\phi\colon M\to N\) is a smooth solution of \eqref{euler-lagrange} and if \(N\)
admits a Killing vector field \(X\), then the following vector field is divergence free
\begin{align}
\label{noether}
J_i:=\delta_1\langle d\phi(e_i),X(\phi)\rangle
+\delta_2\langle\tau(\phi),\bar\nabla_{e_i}X(\phi)\rangle-\delta_2\langle\bar\nabla_{e_i}\tau(\phi),X(\phi)\rangle.
\end{align}
We may rewrite \eqref{noether} as 
\begin{align*}
J_i:=\delta_1\langle d\phi(e_i),X(\phi)\rangle
+\delta_2\nabla_{e_i}\langle\tau(\phi),X(\phi)\rangle-2\delta_2\langle\bar\nabla_{e_i}\tau(\phi),X(\phi)\rangle.
\end{align*}
Recall that for a spherical target the tension field acquires the simple form
\begin{align*}
d\iota(\tau(\phi))=\Delta\varphi+|d\varphi|^2\varphi,
\end{align*}
where \(\iota\colon S^n\to\R^{n+1}\).
Using the expression for the Levi-Civita connection on \(S^n\subset\R^{n+1}\) we find
\begin{align*}
\bar\nabla\tau(\phi)=\nabla\Delta\varphi+\nabla(|d\varphi|^2\varphi)+\langle\Delta\varphi,d\varphi\rangle\varphi.
\end{align*}
This allows us to infer that
\begin{align*}
\langle\bar\nabla\tau(\phi),X(\phi)\rangle=\langle\nabla\Delta\varphi,X(\varphi)\rangle+|d\varphi|^2\langle\nabla\varphi,X(\varphi)\rangle.
\end{align*}
In addition, as \(\varphi\perp X(\varphi)\), we have
\begin{align*}
\langle\tau(\phi),X(\phi)\rangle=\langle\Delta\varphi,X(\varphi)\rangle.
\end{align*}
Combining these equations we find
\begin{align*}
J_i(\varphi)=\delta_2(\langle\nabla_{e_i}\Delta\varphi,X(\varphi)\rangle-\langle\Delta\varphi,\nabla_{e_i} X(\varphi)\rangle+2|d\varphi|^2\langle\nabla_{e_i}\varphi,X(\varphi)\rangle)
-\delta_1\langle\nabla_{e_i}\varphi,X(\varphi)\rangle,
\end{align*}
which is exactly \eqref{conservation-law-sphere-a}.
\end{Bem}

\subsection{Regularity for a spherical target}
In this section we study the regularity of weak solutions of \eqref{euler-lagrange-sphere}
in the case of a spherical target.
Instead of using the Euler-Lagrange equation to obtain some regularity result we will use the 
conserved quantity \eqref{conservation-law-sphere-b} where we follow
the ideas used for intrinsic biharmonic maps from \cite[Theorem A]{MR2094320}.
Moreover, we will apply the regularity theory for extrinsic \cite{MR1692148} and intrinsic biharmonic maps \cite{MR2375314}
and some technical tools to handle lower order terms that have been established in \cite{MR3772035}.

First, let us recall the following 
\begin{Dfn}
For a given open subset \(D\subset\R^m\) and \(1\leq p<\infty\), \(0<\lambda\leq m\)
the Morrey space \(M^{p,\lambda}(D)\) is defined as follows
\begin{align*}
M^{p,\lambda}(D):=\left\{f\in L^p(D)\mid \|f\|^p_{M^{p,\lambda}(D)}:=\sup_{B_r\subset D}\{r^{\lambda-m}\int_{B_r}|f|^p(y)dy\}<\infty\right\}.
\end{align*}
\end{Dfn}
Note that \(M^{p,m}(D)=L^p(D)\).

In the following let \(\Omega\) be an open subset of \(\R^m\) and \(B_r(x)\) the Euclidean ball with radius \(r\) around the point \(x\).
\begin{Dfn}
A map \(\varphi\in W^{2,2}(\Omega,S^n)\) is called a \emph{weak interpolating sesqui-harmonic
map} if it solves \eqref{euler-lagrange-sphere} in the sense of distributions.
\end{Dfn}

First, we will give the following \(\epsilon\)-regularity result.
\begin{Satz}
\label{estimate-sphere}
There exists \(\epsilon_0>0\) such that if 
\(\varphi\in W^{2,2}(\Omega,S^n)\) is a weak interpolating sesqui-harmonic map with \(\delta_2\neq 0\)
then for all \(U\subset B_1(x)\subset\Omega\) with 
\begin{align}
\label{assumption-smoothness}
\|\nabla\varphi\|^4_{M^{4,4}(B_1(x))}+\|\nabla^2\varphi\|^2_{M^{2,4}(B_1(x))}\leq \epsilon_0^2
\end{align}
the following estimate holds
\begin{align}
\label{lp-estimate}
\|\nabla\varphi\|_{L^p(U)}\leq C(\|\nabla\varphi\|_{M^{4,4}(B_1(x))}+1),
\end{align}
where \(4<p<\infty\).
\end{Satz}

In order to prove Theorem \ref{estimate-sphere} we will frequently make use of the Newtonian potential \(I_\alpha\)
which is the operator whose convolution kernel is given by \(|x|^{\alpha-m}\) for \(x\in\R^m\).
We have the following classic estimate for the Newtonian potential in Morrey spaces obtained by Adams \cite{MR0458158}:
\begin{Prop}
\label{adams}
Suppose that \(\alpha>0\), \(0<\lambda\leq m\), \(1<p<\frac{\lambda}{\alpha}\) and let \(f\in M^{p,\lambda}(\R^m)\cap L^p(\R^m)\),
then the following inequality holds
\begin{align}
\label{adams-inequality}
\|I_\alpha(f)\|_{M^{p^\ast,\lambda}(\R^m)}\leq C\|f\|_{M^{p,\lambda}(\R^m)},
\end{align}
where \(p^\ast=\frac{\lambda p}{\lambda-\alpha p}\).
\end{Prop}

We will also need the following auxiliary Lemma:
\begin{Lem}
Suppose that \(h\in C^\infty(B_R)\) with \(0<R<1\) is a solution of \(\Delta^2h=0\).
Then for any \(1<p<\infty\) and \(0<\lambda\leq m\) there exists a constant \(C>0\) such that
for any \(x\in B_{R}\) and \(\theta\in (0,\frac{1}{4})\) the following inequality holds
\begin{align}
\label{estimate-biharmonic-function}
\|\nabla h\|_{M^{p,\lambda}(B_{\theta R}(x))}\leq C\theta\|\nabla h\|_{M^{p,\lambda}(B_{R}(x))}.
\end{align}

\end{Lem}
\begin{proof}
A proof can be found in \cite[p. 229]{MR2094320}, see also \cite[Lemma 4.7]{MR2054520}.
\end{proof}

\begin{proof}[Proof of Theorem \ref{estimate-sphere}]
Let \(G_{\Delta}(x)\) be the Green's operator of \(\Delta\) on \(\R^m\), that is
\begin{align*}
G_{\Delta}(x)&=
\begin{cases}
c_m|x|^{2-m}, &0\neq x\in\R^m, ~m\geq 3 \\
c_2\ln|x|, &0\neq x\in \R^2.
\end{cases}
\end{align*}
Moreover, let \(G_{\Delta^2}(x)\) be the Green's operator of \(\Delta^2\) on \(\R^m\), that is
\begin{align*}
G_{\Delta^2}(x)&=
\begin{cases}
d_m|x|^{4-m}, &0\neq x\in\R^m, ~m\geq 5 \\
d_4\ln|x|, &0\neq x\in \R^4.
\end{cases}
\end{align*}
Here, \(c_m,c_2,d_m,d_4\neq 0\) are constants that only depend on \(m\).

In addition, let \(\tilde\varphi\) be an extension of \(\varphi\) to \(\R^m\) which satisfies
\begin{align*}
\|\nabla\tilde\varphi\|_{M^{p,\lambda}(\R^m)}&\leq C \|\nabla\varphi\|_{M^{p,\lambda}(B_R)}, \\
\|\nabla^2\tilde\varphi\|_{M^{p,\lambda}(\R^m)}&\leq C \|\nabla^2\varphi\|_{M^{p,\lambda}(B_R)},
\end{align*}
where \(0<R<1\). Such an extension can be found by employing a cutoff function which
localizes to the ball \(B_R\).

Moreover, we define the following auxiliary functions \(F_i\colon\R^m\to\R^m,i=1,2\) via
\begin{align*}
F_1(x)&:=-\frac{\delta_1}{\delta_2}\int_{\R^m}G_{\Delta}(x-z)\big(\int_{\R^m}G_{\Delta}(y-z)\nabla(\nabla\tilde\varphi\wedge\tilde\varphi)(y)dy\big)dz, \\
F_2(x)&:=2\int_{\R^m}G_{\Delta^2}(x-y)\nabla\big(\Delta\tilde\varphi\wedge\nabla\tilde\varphi-|d\tilde\varphi|^2\nabla\tilde\varphi\wedge\tilde\varphi\big)(y)dy.
\end{align*}
A direct calculation yields
\begin{align*}
\frac{\partial F_1}{\partial x_i}&=\frac{\delta_1}{\delta_2}\int_{\R^m}\frac{\partial}{\partial x_i}G_{\Delta}(x-z)\big(\int_{\R^m}\frac{\partial}{\partial y_j}G_{\Delta}(y-z)\nabla_j\tilde\varphi\wedge\tilde\varphi)(y)dy\big)dz, \\
\frac{\partial F_2}{\partial x_i}&=-2\int_{\R^m}\frac{\partial^2}{\partial x_i\partial x_j}G_{\Delta_2}(x-y)\big(\Delta\tilde\varphi\wedge\nabla_j\tilde\varphi-|d\tilde\varphi|^2\nabla_j\tilde\varphi\wedge\tilde\varphi\big)(y)dy.
\end{align*}

Therefore, we obtain the following inequalities
\begin{align*}
|\nabla F_1|(x)\leq& C\int_{\R^m}|x-z|^{1-m}\big(\int_{\R^m}|y-z|^{1-m}|\nabla\tilde\varphi|(y)dy\big)dz \\
\leq& CI_1(I_1(|\nabla\tilde\varphi|))(x)
\end{align*}
and
\begin{align*}
|\nabla F_2|(x)\leq& C\int_{\R^m}|x-y|^{2-m}(|\Delta\tilde\varphi||\nabla\tilde\varphi|+|\nabla\tilde\varphi|^3)(y)dy \\
\leq & C\big(I_2(|\Delta\tilde\varphi||\nabla\tilde\varphi|)(x)+I_2(|\nabla\tilde\varphi|^3)(x)\big).
\end{align*}

Applying Proposition \ref{adams}, choosing \(\lambda=4, p^\ast=4\) and \(\alpha=1\), we find
\begin{align*}
\|\nabla F_1\|_{M^{4,4}(\R^m)}\leq C\|I_1(I_1(|\nabla\tilde\varphi|))\|_{M^{4,4}(\R^m)}\leq C\|I_1(|\nabla\tilde\varphi|)\|_{M^{2,4}(\R^m)}\leq C\|\nabla\tilde\varphi\|_{M^{\frac{4}{3},4}(\R^m)}.
\end{align*}

Applying Proposition \ref{adams} again we find by choosing \(p=\frac{4}{3}, \alpha=2\) and \(\lambda=4\) that
\begin{align*}
\|\nabla F_2\|_{M^{4,4}(\R^m)}\leq C(\||\Delta\tilde\varphi||\nabla\tilde\varphi|\|_{M^{\frac{4}{3},4}(\R^m)}+\||\nabla\tilde\varphi|^3\|_{M^{\frac{4}{3},4}(\R^m)}).
\end{align*}

Using Hölder's inequality we get the estimates
\begin{align*}
\||\Delta\tilde\varphi||\nabla\tilde\varphi|\|_{M^{\frac{4}{3},4}(\R^m)}&\leq \|\Delta\tilde\varphi\|_{M^{2,4}(\R^m)}\|\nabla\tilde\varphi\|_{M^{4,4}(\R^m)},\\
\||\nabla\tilde\varphi|^3\|_{M^{\frac{4}{3},4}(\R^m)}&\leq \|\nabla\tilde\varphi\|^3_{M^{4,4}(\R^m)}
\end{align*}
and also
\begin{align*}
\|\nabla\tilde\varphi\|_{M^{\frac{4}{3},4}(\R^m)}&\leq C\|\nabla\varphi\|_{M^{\frac{4}{3},4}(B_R)}\\
&\leq C\|\nabla\varphi\|_{M^{4,4}(B_R)}\vol(B_R)^\frac{1}{2}R^{2-\frac{m}{2}}\\
&\leq C\epsilon_0^\frac{1}{2}R^2.
\end{align*}

Setting \(F=F_1+F_2\) we find
\begin{align}
\label{estimate-nabla-F}
\|\nabla F\|_{M^{4,4}(\R^m)}&\leq C(\|\Delta\tilde\varphi\|_{M^{2,4}(\R^m)}\|\nabla\tilde\varphi\|_{M^{4,4}(\R^m)}+\|\nabla\tilde\varphi\|^3_{M^{4,4}(\R^m)}+\|\nabla\tilde\varphi\|_{M^{\frac{4}{3},4}(\R^m)}) \\
\nonumber &\leq C\epsilon_0\|\nabla\varphi\|_{M^{4,4}(B_R)}+C\|\nabla\varphi\|_{M^{\frac{4}{3},4}(B_R)}\\
\nonumber &\leq C\epsilon_0\|\nabla\varphi\|_{M^{4,4}(B_R)}+CR^2,
\end{align}
where we have used the assumption \eqref{assumption-smoothness}.

Now, we consider the following distributional version of the Hodge decomposition of the one-form \(\nabla\tilde\varphi\wedge\tilde\varphi\)
which was given in \cite[Theorem 6.1]{MR1208562}:
There exist \(\Phi\in W^{1,4}(\R^m)\) and \(\psi\in W^{1,4}(\R^m,\Lambda^2(\R^m))\) such that
\begin{align}
\label{hodge-decomposition}
\nabla\tilde\varphi\wedge\tilde\varphi=d\Phi+d^\ast\psi,\qquad d\psi=0
\end{align}
and also
\begin{align*}
\|\nabla\Phi\|_{L^4(\R^m)}+\|\nabla\psi\|_{L^4(\R^m)}\leq C\|\nabla\tilde\varphi\|_{L^4(\R^m)}.
\end{align*}
Here, \(d^\ast\) denotes the codifferential.
Applying \(\Delta d^\ast\) to \eqref{hodge-decomposition} we obtain
\begin{align*}
\Delta^2\Phi=\nabla\big(2\Delta\varphi\wedge\nabla\varphi-2|d\varphi|^2\nabla\varphi\wedge\varphi-\frac{\delta_1}{\delta_2}\nabla\varphi\wedge\varphi\big)\qquad \text{ in } B_R,
\end{align*}
where we made use of \eqref{conservation-law-sphere-b}.
Moreover, we obtain
\begin{align*}
\Delta^2F=\nabla\big(2\Delta\tilde\varphi\wedge\nabla\tilde\varphi-2|d\tilde\varphi|^2\nabla\tilde\varphi\wedge\tilde\varphi-\frac{\delta_1}{\delta_2}\nabla\tilde\varphi\wedge\tilde\varphi\big)\qquad \text{ in } \R^m.
\end{align*}
Consequently, we find
\begin{align*}
\Delta^2(\Phi-F)=0\qquad \text{ in } B_R.
\end{align*}
From \eqref{estimate-biharmonic-function} we may then deduce that for any \(\theta\in (0,\frac{1}{4})\) we have
\begin{align*}
\|\nabla(\Phi-F)\|_{M^{4,4}(B_{\theta R})}\leq C\theta \|\nabla(\Phi-F)\|_{M^{4,4}(B_{R})}
\end{align*}
and together with \eqref{estimate-nabla-F} this yields
\begin{align}
\label{estimate-bigphi}
\|\nabla\Phi\|_{M^{4,4}(B_{\theta R})}\leq C(\theta+\epsilon_0)\|\nabla\varphi\|_{M^{4,4}(B_{R})}+CR^2.
\end{align}
Applying the exterior derivative \(d\) on both sides of \eqref{hodge-decomposition} we get the following
equation for \(\psi\) 
\begin{align*}
\Delta\psi=d\tilde\varphi\wedge d\tilde\varphi.
\end{align*}
Using the explicit formula for the solution of the Poisson equation in \(\R^m\) we find the following
estimate
\begin{align*}
|\nabla\psi|(x)\leq CI_1(|\nabla\tilde\varphi|^2)(x).
\end{align*}
Thanks to \eqref{adams-inequality} we obtain
\begin{align}
\label{estimate-psi}
\|\nabla\psi\|_{M^{4,4}(\R^m)}\leq C\||\nabla\tilde\varphi|^2\|_{M^{2,4}(\R^m)}\leq C\epsilon_0\|\nabla\varphi\|_{M^{4,4}(B_R)}.
\end{align}
Combining \eqref{hodge-decomposition} with \eqref{estimate-bigphi}, \eqref{estimate-psi} we can conclude that
\begin{align*}
\|\nabla\varphi\wedge\varphi\|_{M^{4,4}(B_{\theta R})}\leq C(\epsilon_0+\theta)\|\nabla\varphi\|_{M^{4,4}(B_{R})}+CR^2.
\end{align*}
As we are considering a spherical target we have 
\begin{align*}
|\nabla\varphi\wedge\varphi|^2=|\nabla\varphi|^2.
\end{align*}

Now, fix any \(\beta\in(0,1)\), then we can find \(\theta\in (0,\frac{1}{4})\) such that \(C\theta\leq\frac{1}{2}\theta^\beta\).
Now, we choose \(\epsilon_0\) small enough such that \(C\epsilon_0\leq\frac{1}{2}\theta^\beta\).
This leads us to the following inequality
\begin{align}
\label{iteration-start}
\|\nabla\varphi\|_{M^{4,4}(B_{\theta R})}=\|\nabla\varphi\wedge\varphi\|_{M^{4,4}(B_{\theta R})}
\leq \theta^\beta\|\nabla\varphi\|_{M^{4,4}(B_R)}+CR^2.
\end{align}

At this point we are starting an iteration procedure as in \cite[p. 194]{MR3772035}.
To this end we consider \(R<1-\delta\) for some small \(\delta>0\).
Then, for any \(0<r<R\), there exists \(k\in\N\) such that \(\theta^{k+1}R<r\leq\theta^kR\).
Iterating \eqref{iteration-start} we obtain
\begin{align*}
\|\nabla\varphi\|_{M^{4,4}(B_{r})}&\leq\|\nabla\varphi\|_{M^{4,4}(B_{\theta^kR})} \\
&\leq \theta^\beta\|\nabla\varphi\|_{M^{4,4}(B_{\theta^{k-1}R})}+C(R\theta^{k-1})^2\\
&\leq \theta^{2\beta}\|\nabla\varphi\|_{M^{4,4}(B_{\theta^{k-2}R})}+C(R\theta^{k-1})^2+C\theta^\beta(R\theta^{k-2})^2\\
&\leq\ldots\\
&\leq \theta^{k\beta}\|\nabla\varphi\|_{M^{4,4}(B_R)}+CR^2\theta^{(k-1)\beta}\sum_{j=0}^{k-1}\theta^{(2-\beta)j}\\
&\leq \frac{1}{\theta^\beta}\big(\frac{r}{R}\big)^\beta\|\nabla\varphi\|_{M^{4,4}(B_R)}
+CR^2\theta^{(k-1)\beta}\frac{1-\theta^{(2-\beta)k}}{1-\theta^{2-\beta}}\\
&= \frac{1}{\theta^\beta}\big(\frac{r}{R}\big)^\beta\|\nabla\varphi\|_{M^{4,4}(B_R)}
+C\frac{R^{2-\beta}}{\theta^{2\beta}}
\frac{1-\theta^{(2-\beta)k}}{1-\theta^{2-\beta}}(\theta^{k+1}R)^\beta\\
&\leq \frac{1}{\theta^\beta}\big(\frac{r}{R}\big)^\beta\|\nabla\varphi\|_{M^{4,4}(B_R)}+Cr^\beta,
\end{align*}
where we used that \(R<1\) in the last step.

For \(\delta<\frac{1}{4}\) and \(R=\frac{1}{2}\) we may conclude that
\begin{align*}
\|\nabla\varphi\|_{M^{4,4-4\beta}(B_{\frac{1}{4}})}&\leq C\|\nabla\varphi\|_{M^{4,4}(B_1)}+C
\end{align*}
and thus \(\nabla\varphi\in M^{4,4-4\beta}(B_{\frac{1}{4}})\).

From now on we will assume that \(0<R<\frac{1}{8}\).
As a next step we want to improve the integrability of \(\nabla\varphi\). To this end we 
apply Proposition \ref{adams}, choosing \(\lambda=4-\frac{4}{3}\beta,p^\ast=\frac{4-\frac{4}{3}\beta}{1-\beta}\) and \(\alpha=1\), and find
\begin{align*}
\|\nabla F_1\|_{M^{\frac{4-\frac{4}{3}\beta}{1-\beta},4-\frac{4}{3}\beta}(\R^m)}
&\leq C\|I_1(I_1(|\nabla\tilde\varphi|))\|_{M^{\frac{4-\frac{4}{3}\beta}{1-\beta},4-\frac{4}{3}\beta}(\R^m)} \\
&\leq C\|I_1(|\nabla\tilde\varphi|)\|_{M^{\frac{4-\frac{4}{3}\beta}{2-\beta},4-\frac{4}{3}\beta}(\R^m)} \\
&\leq C\|\nabla\tilde\varphi\|_{M^{\frac{4}{3},4-\frac{4}{3}\beta}(\R^m)}\\
&\leq C\|\nabla\varphi\|_{M^{\frac{4}{3},4-\frac{4}{3}\beta}(B_R)}\\
&\leq C\epsilon_0^\frac{1}{2}R^{2-\beta}.
\end{align*}

In addition, applying Proposition \ref{adams} with \(p=\frac{4}{3},\alpha=2\) and \(\lambda=4-\frac{4}{3}\beta\) once more, we find
\begin{align*}
\|\nabla F_2\|_{M^{\frac{4-\frac{4}{3}\beta}{1-\beta},4-\frac{4}{3}\beta}(\R^m)}&\leq C(\||\Delta\tilde\varphi||\nabla\tilde\varphi|\|_{M^{\frac{4}{3},4-\frac{4}{3}\beta}(\R^m)}
+\||\nabla\tilde\varphi|^3\|_{M^{\frac{4}{3},4-\frac{4}{3}\beta}(\R^m)})\\
&\leq C(\|\Delta\tilde\varphi\|_{M^{2,4}(\R^m)}\|\nabla\tilde\varphi\|_{M^{4,4-4\beta}(\R^m)}
+\|\nabla\tilde\varphi\|^2_{M^{4,4}(\R^m)}\|\nabla\tilde\varphi\|_{M^{4,4-4\beta}(\R^m)})\\
&\leq C(\|\Delta\varphi\|_{M^{2,4}(B_R)}\|\nabla\varphi\|_{M^{4,4-4\beta}(B_R)}
+\|\nabla\varphi\|^2_{M^{4,4}(B_R)}\|\nabla\varphi\|_{M^{4,4-4\beta}(B_R)})\\
&\leq C\epsilon_0\|\nabla\varphi\|_{M^{4,4-4\beta}(B_R)}.
\end{align*}

Combining both inequalities we find
\begin{align}
\label{estimate-nabla-F-b}
\|\nabla F\|_{M^{\frac{4-\frac{4}{3}\beta}{1-\beta},4-\frac{4}{3}\beta}(\R^m)}\leq C\epsilon_0\|\nabla\varphi\|_{M^{4,4-4\beta}(B_R)}+CR^{2-\beta}.
\end{align}

Recall that
\begin{align*}
\Delta^2(\Phi-F)=0\qquad \text{ in } B_R.
\end{align*}
Using \eqref{estimate-biharmonic-function} together with \eqref{estimate-nabla-F-b} this yields
\begin{align*}
\|\nabla\Phi\|_{M^{\frac{4-\frac{4}{3}\beta}{1-\beta},4-\frac{4}{3}\beta}(B_{\theta R})}\leq C(\theta+\epsilon_0)\|\nabla\varphi\|_{M^{4,4-4\beta}(B_R)}+CR^2.
\end{align*}

Applying \(d\) first and \(\nabla\) afterwards in \eqref{hodge-decomposition} we can derive the following estimate
\begin{align*}
|\nabla\psi|(x)\leq C I_2(|\nabla^2\tilde\varphi||\nabla\tilde\varphi|)(x).
\end{align*}

As \(\beta\in(0,1)\) this allows us to calculate
\begin{align*}
\|\nabla\psi\|_{M^{\frac{4-\frac{4}{3}\beta}{1-\beta},4-\frac{4}{3}\beta}(\R^m)}
&\leq C\|I_2(|\nabla^2\tilde\varphi||\nabla\tilde\varphi|)\|_{M^{\frac{4-\frac{4}{3}\beta}{1-\beta},4-\frac{4}{3}\beta}(\R^m)} \\
&\leq C\||\nabla^2\tilde\varphi||\nabla\tilde\varphi|\|_{M^{\frac{4}{3},4-\frac{4}{3}\beta}(\R^m)}\\
&\leq C\|\nabla^2\varphi\|_{M^{2,4}(B_R)}\|\nabla\varphi\|_{M^{4,4-4\beta}(B_R)}.
\end{align*}

We may conclude that
\begin{align}
\label{morrey-estimate}
\|\nabla\varphi\|_{M^{\frac{4-\frac{4}{3}\beta}{1-\beta},4-\frac{4}{3}\beta}(B_\frac{1}{32})}
&=\|\nabla\varphi\wedge\varphi\|_{M^{\frac{4-\frac{4}{3}\beta}{1-\beta},4-\frac{4}{3}\beta}(B_\frac{1}{32})}\\
\nonumber &\leq\|\nabla\psi\|_{M^{\frac{4-\frac{4}{3}\beta}{1-\beta},4-\frac{4}{3}\beta}(B_\frac{1}{32})}
+\|\nabla\Phi\|_{M^{\frac{4-\frac{4}{3}\beta}{1-\beta},4-\frac{4}{3}\beta}(B_\frac{1}{32})} \\
&\nonumber\leq C(\|\nabla\varphi\|_{M^{4,4}(B_1)}+1).
\end{align}

Note that \eqref{morrey-estimate} holds for any \(0<\beta<1\).
Moreover, as
\begin{align*}
\lim_{\beta\to 1}\frac{4-\frac{4}{3}\beta}{1-\beta}=\infty
\end{align*}
and also
\begin{align*}
M^{\frac{4-\frac{4}{3}\beta}{1-\beta},4-\frac{4}{3}\beta}(B_\frac{1}{32})\hookrightarrow L^{\frac{4-\frac{4}{3}\beta}{1-\beta}}(B_\frac{1}{32})
\end{align*}
we can conclude that \(\nabla\varphi\in L^p(B_\frac{1}{32})\) for any \(4<p<\infty\)
completing the proof.
\end{proof}

Using the regularity theory for biharmonic maps from four-dimensional domains \cite{MR2375314} we 
can now give the following regularity result:
\begin{Satz}
Let \(\varphi\colon\Omega\to S^n\) be a weak solution of \eqref{euler-lagrange-sphere} with \(\delta_2\neq 0\)
that satisfies \eqref{assumption-smoothness} and suppose that \(\dim M=4\).
Then \(\varphi\) is smooth, that is \(\varphi\in C^\infty(\Omega,S^n)\).
\end{Satz}
\begin{proof}
Thanks to the estimate \eqref{lp-estimate} the Morrey Lemma yields that \(\varphi\in C^\alpha(B_{\frac{1}{2}}(x),S^{n})\) for some \(\alpha>0\).
The result now follows from the regularity theory for biharmonic maps to spheres \cite[Theorem 5.1]{MR1692148}, \cite[Theorem 3.1]{MR2375314},
see also \cite{MR1809291,MR2094552}.
\end{proof}

\begin{Bem}
If \(\dim M=4\) the smallness condition \eqref{assumption-smoothness} reads
\begin{align*}
\int_{B_1(x)}(|\nabla^2\varphi|^2+|\nabla\varphi|^4)d\mu\leq\epsilon_0^2.
\end{align*}
\end{Bem}

\begin{Bem}
Using the refined techniques for biharmonic maps from a four-dimensional domain to an arbitrary Riemannian manifold developed in \cite{MR2398228,MR2054520}
together with the results obtained in this article it should be possible to also establish 
the regularity of weak solutions for interpolating sesqui-harmonic maps to an arbitrary target manifold.
\end{Bem}

\subsection{Some remarks on interpolating sesqui-harmonic immersions to spheres}
In this subsection we present a classification result for interpolating sesqui-harmonic
immersions to spaces of positive constant sectional curvature \(K\), where we apply ideas that
have been used to study triharmonic immersions from \cite[Section 3]{MR2869168}.

Recall that if \(\phi\) is an isometric immersion then \(\langle d\phi(X),\tau(\phi)\rangle=0\) for 
all \(X\in\Gamma(TM)\).

Due to our assumptions on the map \(\phi\) and the geometry of the target we find
\begin{align*}
R^N(d\phi(e_i),\tau(\phi))d\phi(e_i)=K(\langle d\phi(e_i),\tau(\phi)\rangle d\phi(e_i)-|d\phi|^2\tau(\phi))=-mK\tau(\phi),
\end{align*}
where \(m=\dim M\).

Hence, under the above assumptions, the Euler-Lagrange equation \eqref{euler-lagrange} acquires the simple form 
\begin{align}
\label{isometric-immersion-sphere}
\bar\Delta\tau(\phi)=(-mK+\frac{\delta_1}{\delta_2})\tau(\phi).
\end{align}

For solutions of \eqref{isometric-immersion-sphere} we will prove the following result:

\begin{Satz}
\label{theorem-immersion-sphere}
Let \(\phi\colon M\to S^n\) be a smooth solution of \eqref{isometric-immersion-sphere}
and \(-mK+\frac{\delta_1}{\delta_2}\geq 0\), where \(m=\dim M\).
Then the following statements hold
\begin{enumerate}
 \item If \(M\) is closed then \(\phi\) must be harmonic.
 \item If \(M\) is complete, non-compact and \(\phi\) has finite bienergy, that is \(E_2(\phi)<\infty\), then
  \(\phi\) must be harmonic.
\end{enumerate}
Here, \(K\) represents the curvature of \(S^n\).
\end{Satz}
\begin{proof}
To prove the first assertion we test \eqref{isometric-immersion-sphere} with \(\tau(\phi)\).
After using integration by parts we obtain
\begin{align*}
-\int_M|\bar\nabla\tau(\phi)|^2\dv=(-mK+\frac{\delta_1}{\delta_2})\int_M|\tau(\phi)|^2\dv\geq 0
\end{align*}
which yields \(\bar\nabla\tau(\phi)=0\).
As \(\phi\) is an isometric immersion we can conclude that
\begin{align*}
0=\langle d\phi,\bar\nabla\tau(\phi)\rangle=\nabla\underbrace{\langle d\phi,\tau(\phi)\rangle}_{=0}-|\tau(\phi)|^2
\end{align*}
finishing the proof.

In order to prove the second claim
we choose a cutoff function  \(0\leq\eta\leq 1\) on \(M\) that satisfies
\begin{align*}
\eta(x)=1\textrm{ for } x\in B_R(x_0),\qquad \eta(x)=0\textrm{ for } x\in B_{2R}(x_0),\qquad |\nabla\eta|\leq\frac{C}{R}\textrm{ for } x\in M,
\end{align*}
where \(B_R(x_0)\) denotes the geodesic ball around the point \(x_0\) with radius \(R\).
Testing \eqref{isometric-immersion-sphere} with \(\tau(\phi)\eta^2\) and applying the assumptions we obtain
\begin{align*}
0\leq\int_M\langle\bar\Delta\tau(\phi),\tau(\phi)\rangle\eta^2\dv
=-\int_M|\bar\nabla\tau(\phi)|^2\eta^2\dv-2\int_M\langle\bar\nabla\tau(\phi),\tau(\phi)\rangle\eta\nabla\eta\dv.
\end{align*}
As the bienergy of \(\phi\) is finite by assumption we may conclude that
\begin{align*}
\int_M|\bar\nabla\tau(\phi)|^2\dv\leq \frac{C}{R^2}E_2(\phi)\to 0 \text{ as } R\to\infty
\end{align*}
and the claim follows by the same method as in the first case.
\end{proof}

\begin{Bem}
By inspecting the assumptions of Theorem \ref{theorem-immersion-sphere} we realize
that we have to require that \(mK\leq\delta_1/\delta_2\). We realize that \(\delta_1\)
and \(\delta_2\) need to have the same sign and that their ratio needs to be sufficiently large.
Moreover, it is necessary that \(\delta_1\neq 0\).
\end{Bem}

\section{A growth formula and a classification result}
In this section we prove a classification result for interpolating sesqui-harmonic maps
from \(\R^m\) under certain boundedness assumptions. In addition, we also establish a growth formula 
for interpolating sesqui-harmonic maps from \(\R^m\). Both results make use of the stress-energy tensor
associated to interpolating sesqui-harmonic maps which is obtained by varying the functional \eqref{energy-functional}
with respect to the metric on the domain.
This tensor was derived in \cite[Proposition 2.6]{MR4058514} where it was also shown that
it is divergence free whenever we have a solution of \eqref{euler-lagrange}.

In terms of a local orthonormal basis \(\{e_i\},i=1,\ldots,m\) we may express the stress-energy tensor as
\begin{align}
\label{energy-momentum-coordinates}
S_{ij}=&\delta_1(2\langle d\phi(e_i),d\phi(e_j)\rangle-g_{ij}|d\phi|^2) \\
\nonumber &+\delta_2(-|\tau(\phi)|^2+2\nabla_{e_k}\langle d\phi(e_k),\tau(\phi)\rangle)g_{ij}
-2\delta_2(\langle d\phi(e_i),\bar\nabla_{e_j}\tau(\phi)\rangle+\langle d\phi(e_j),\bar\nabla_{e_i}\tau(\phi)\rangle).
\end{align}
The stress-energy tensor for polyharmonic maps was recently studied in \cite{MR4007262}.

\subsection{A classification result}
The classification result that we will derive in the following does not necessarily require that
we are considering a smooth solution of \eqref{euler-lagrange}. In order to state the result we 
give the following 
\begin{Dfn}
A solution \(\phi\in W^{2,2}_{loc}(M,N)\cap L^\infty_{loc}(M,N)\) of \eqref{euler-lagrange} is called 
\emph{stationary} if it is also a critical point of \eqref{energy-functional} with respect to variations
of the metric on the domain, that is
\begin{align}
0=\int_Mk^{ij}S_{ij}\dv,
\end{align}
where \(k_{ij}\) is a smooth symmetric 2-tensor.
\end{Dfn}

We obtain the following vanishing result for stationary solutions of \eqref{euler-lagrange}:
\begin{Satz}
\label{theorem-liouville}
Let \(\phi\in W^{2,2}_{loc}(\R^m,N)\cap L^\infty_{loc}(\R^m,N)\) be a stationary solution of \eqref{euler-lagrange} with \(\delta_1,\delta_2\neq 0\).
Moreover, suppose that
\begin{align}
\label{finite-energy}
\int_{\R^m}(|d\phi|^2+|\bar\nabla d\phi|^2)\dv<\infty.
\end{align}
Then we obtain the following kind of classification result:
\begin{enumerate}
 \item If \(m=1\) and \(\delta_1,\delta_2>0\) then \(\phi\) is trivial.
 \item If \(m=2\) and \(\delta_2>0\) then \(\phi\) is harmonic.
 \item If \(m=3\) and \(\delta_1<0,\delta_2>0\) then \(\phi\) is trivial.
 \item If \(m=4\) and \(\delta_1<0\) then \(\phi\) is trivial.
 \item If \(m\geq 5\) and \(\delta_1,\delta_2<0\) then \(\phi\) is trivial.
\end{enumerate}

\end{Satz}
\begin{proof}
Let \(\eta\in C_0^\infty(\mathbb{R})\) be a smooth cutoff function satisfying \(\eta=1\) for \(r\leq R\),
\(\eta=0\) for \(r\geq 2R\) and \(|\eta^{(i)}(r)|\leq\frac{C}{R^i},i=1,2\). In addition, 
we choose \(Y(x):=x\eta(r)\in C_0^\infty(\mathbb{R}^m,\mathbb{R}^m)\) with \(r=|x|\).
Hence, we find
\[
k_{ij}=\frac{\partial Y_i}{\partial x^j}=\delta_{ij}\eta(r)+\frac{x_i x_j}{r}\eta'(r).
\]
By assumption the map \(\phi\) is stationary which means that 
\begin{align*}
0=\int_Mk^{ij}S_{ij}\dv.
\end{align*}
Now, a direct computation yields
\begin{align*}
\int_{\R^m}\delta_{ij}\eta(r)S_{ij}\dv=&\delta_1(2-m)\int_{\R^m}\eta(r)|d\phi|^2\dv+\delta_2(4-m)\int_{\R^m}\eta(r)|\tau(\phi)|^2\dv \\
&-\delta_2(2m-4)\int_{\R^m}\eta'(r)\frac{x_i}{r}\langle d\phi(e_i),\tau(\phi)\rangle\dv.
\end{align*}
Moreover, we find
\begin{align*}
\int_{\R^m}\eta'(r)\frac{x_ix_j}{r}S_{ij}\dv=&\delta_1\int_{\R^m}\eta'(r)r\big(2|d\phi(\partial_r)|^2-|d\phi|^2\big)\dv\\
&+\delta_2\int_{\R^m}\eta'(r)r(-|\tau(\phi)|^2+2\nabla_{e_k}\langle d\phi(e_k),\tau(\phi)\rangle)\dv\\
&-4\delta_2\int_{\R^m}\eta'(r)\frac{x_ix_j}{r}\langle d\phi(e_i),\bar\nabla_{e_j}\tau(\phi)\rangle\dv.
\end{align*}
By integration by parts we obtain
\begin{align*}
\int_{\R^m}\eta'(r)r\nabla_{e_k}\langle d\phi(e_k),\tau(\phi)\rangle\dv=&
-\int_{\R^m}\big(\eta''(r)x_k+\eta'(r)\frac{x_k}{r}\big)\langle d\phi(e_k),\tau(\phi)\rangle\dv
\end{align*}
and also
\begin{align*}
\int_{\R^m}\eta'(r)\frac{x_ix_j}{r}\langle d\phi(e_i),\bar\nabla_{e_j}\tau(\phi)\rangle\dv=&
-\int_{\R^m}\big(\eta''(r)x_i+2\eta'(r)\frac{x_i}{r}\big)\langle d\phi(e_i),\tau(\phi)\rangle\dv \\
&-\int_{\R^m}\eta'(r)\frac{x_ix_j}{r}\langle\bar\nabla_{e_j} d\phi(e_i),\tau(\phi)\rangle\dv.
\end{align*}

This leads us to the following equality
\begin{align}
\label{vanishing-a}
\delta_1(2-m)&\int_{\R^m}\eta(r)|d\phi|^2\dv+\delta_2(4-m)\int_{\R^m}\eta(r)|\tau(\phi)|^2\dv\\
\nonumber=&\delta_1\int_{\R^m}\eta'(r)r\big(|d\phi|^2-2|d\phi(\partial_r)|^2\big)\dv
+\delta_2(2m-10)\int_{\R^m}\eta'(r)\frac{x_i}{r}\langle d\phi(e_i),\tau(\phi)\rangle\dv \\
\nonumber&+\delta_2\int_{\R^m}\eta'(r)r|\tau(\phi)|^2\dv
-2\delta_2\int_{\R^m}\eta''(r)x_i\langle d\phi(e_i),\tau(\phi)\rangle\dv \\
\nonumber&-4\delta_2\int_{\R^m}\eta'(r)\frac{x_ix_j}{r}\langle\bar\nabla_{e_j} d\phi(e_i),\tau(\phi)\rangle\dv.
\end{align}
We can control the right-hand side of \eqref{vanishing-a} as follows
\begin{align*}
\int_{\R^m} r\eta'(r)(|d\phi|^2-2|d\phi(\partial_r)|^2)\dv
&\leq C\int_{B_{2R}\setminus B_R}|d\phi|^2d\mu\rightarrow 0 \text{  as  } R\to\infty,\\
\int_{\R^m}\eta'(r)\frac{x_i}{r}\langle d\phi(e_i),\tau(\phi)\rangle\dv
&\leq \frac{C}{R}\int_{\R^m}(|d\phi|^2+|\tau(\phi)|^2)\dv\rightarrow 0 \text{  as  } R\to\infty,\\
\int_{\R^m}\eta'(r)r|\tau(\phi)|^2\dv
&\leq C\int_{B_{2R}\setminus B_R}|\tau(\phi)|^2d\mu\rightarrow 0 \text{  as  } R\to\infty,\\
\int_{\R^m}\eta''(r)x_i\langle d\phi(e_i),\tau(\phi)\rangle\dv
&\leq \frac{C}{R}\int_{\R^m}(|d\phi|^2+|\tau(\phi)|^2)\dv\rightarrow 0 \text{  as  } R\to\infty, \\
\int_{\R^m}\eta'(r)\frac{x_ix_j}{r}\langle\nabla_{e_j} d\phi(e_i),\tau(\phi)\rangle\dv
&\leq C\int_{B_{2R}\setminus B_R}|\bar\nabla d\phi|^2d\mu\rightarrow 0 \text{  as  } R\to\infty.
\end{align*}
Here, we used the finiteness assumption \eqref{finite-energy} and that \(|\tau(\phi)|^2\leq m|\bar\nabla d\phi|^2\).
Consequently, we obtain from \eqref{vanishing-a} after taking the limit \(R\to\infty\) that
\begin{align*}
\delta_1(2-m)&\int_{\R^m}|d\phi|^2\dv+\delta_2(4-m)\int_{\R^m}|\tau(\phi)|^2\dv\leq 0.
\end{align*}
The result follows from this formula by performing a case by case analysis.
\end{proof}

Note that Theorem \ref{theorem-liouville} generalizes some vanishing results for biharmonic maps obtained in \cite{MR2604617}.

\subsection{Monotonicity formulas}
A monotonicity formula for biharmonic immersions satisfying an additional assumption was established in \cite[Theorem 5.1]{MR3178159}.
Without the assumption that \(\phi\) is an immersion it is rather cumbersome to derive a monotonicity formula.

In the following we will derive a monotonicity formula for solutions of \eqref{euler-lagrange}
where, for simplicity, we will stick to the case of a Euclidean domain.

Let us recall the following facts:
A vector field \(X\) is called \emph{conformal} if
\[
\mathcal{L}_Xg=fg,
\]
where \(\mathcal{L}_X\) denotes the Lie derivative of the metric with respect to \(X\) and
\(f\colon M\to\mathbb{R}\) is a smooth function.
\begin{Lem}
Let \(T\) be a symmetric 2-tensor. For any vector field \(X\) the following formula holds
\begin{equation*}
\operatorname{div}(\iota_X T)=\iota_X\operatorname{div} T+\langle T,\nabla X\rangle.
\end{equation*}
If \(X\) is a conformal vector field, then the second term on the right hand side acquires the form
\begin{equation}
\label{conformal-vf}
\langle T,\nabla X\rangle=\frac{1}{m}\operatorname{div}X\tr T.
\end{equation}
\end{Lem}
By integrating over a compact region \(U\), making use of Stokes theorem, we obtain
\begin{Lem}
Let \((M,g)\) be a Riemannian manifold and \(U\subset M\) be a compact region with smooth boundary.
Then, for any symmetric \(2\)-tensor and any vector field \(X\) the following formula holds
\begin{equation}
\label{gauss-tensor-formula}
\int_{\partial U}T(X,\nu)d\sigma=\int_U\iota_X\operatorname{div} Td\mu+\int_U\langle T,\nabla X\rangle d\mu,
\end{equation}
where \(\nu\) denotes the normal to \(U\).
The same formula holds for a conformal vector field \(X\) if we replace the second 
term on the right hand by \eqref{conformal-vf}.
\end{Lem}

\begin{Lem}
Let \(\phi\colon\R^m\to N\) be a smooth solution of \eqref{euler-lagrange}.
Then the following equality holds
\begin{align}
\label{identity-general-rm}
r^{m-2}\frac{d}{dr}r^{2-m}&\delta_1\int_{B_r}|d\phi|^2d\mu+r^{m-4}\frac{d}{dr}r^{4-m}\delta_2\int_{B_r}|\tau(\phi)|^2d\mu \\
\nonumber=&2\delta_1\int_{\partial B_r}|d\phi(\partial_r)|^2d\sigma
+\frac{2\delta_2(2-m)}{r}\int_{\partial B_r}\langle d\phi(\partial_r),\tau(\phi)\rangle d\sigma \\
&+2\delta_2\int_{\partial B_r}\nabla_{e_k}\langle d\phi(e_k),\tau(\phi)\rangle d\sigma
\nonumber-4\delta_2\int_{\partial B_r}\langle d\phi(\partial_r),\bar\nabla_{\partial_r}\tau(\phi)\rangle d\sigma,
\end{align}
where \(B_r\) denotes the geodesic ball of radius \(r\) in \(\R^m\).
\end{Lem}
\begin{proof}
We choose the conformal vector field \(X=r\frac{\partial}{\partial r}\) which satisfies \(\operatorname{div} X=m\) and apply \eqref{gauss-tensor-formula} to
the stress-energy tensor \eqref{energy-momentum-coordinates}.
We obtain the following equality
\begin{align*}
\delta_1(2-m)&\int_{B_r}|d\phi|^2d\mu+\delta_2m\int_{B_r}2(\nabla\langle d\phi,\tau(\phi)\rangle-|\tau(\phi)|^2)d\mu
-4\delta_2\int_{B_r}\langle d\phi,\bar\nabla\tau(\phi)\rangle d\mu \\
=&\delta_1r\int_{\partial B_r}(2|d\phi(\partial_r)|^2-|d\phi|^2)d\sigma \\
&+\delta_2r\int_{\partial B_r}(-|\tau(\phi)|^2+2\nabla\langle d\phi,\tau(\phi)\rangle-4\langle d\phi(\partial_r),\bar\nabla_{\partial_r}\tau(\phi)\rangle)d\sigma.
\end{align*}
As a next step we employ integration by parts to deduce
\begin{align*}
\int_{B_r}\langle d\phi,\bar\nabla\tau(\phi)\rangle d\mu=&-\int_{B_r}|\tau(\phi)|^2d\mu+\int_{\partial B_r}\langle\tau(\phi),d\phi(\partial_r)\rangle d\sigma, \\
\int_{B_r}\nabla\langle d\phi,\tau(\phi)\rangle d\mu=&\int_{\partial B_r}\langle\tau(\phi),d\phi(\partial_r)\rangle d\sigma.
\end{align*}
This leads us to
\begin{align*}
\delta_1(2-m)&\int_{B_r}|d\phi|^2d\mu+\delta_1r\int_{\partial B_r}|d\phi|^2d\sigma+\delta_2(4-m)\int_{B_r}|\tau(\phi)|^2d\mu+\delta_2r\int_{\partial B_r}|\tau(\phi)|^2d\sigma \\
=&2\delta_1r\int_{\partial B_r}|d\phi(\partial_r)|^2d\sigma 
+2\delta_2(2-m)\int_{\partial B_r}\langle d\phi(\partial_r),\tau(\phi)\rangle d\sigma \\
&+\delta_2r\int_{\partial B_r}\big(2\nabla\langle d\phi,\tau(\phi)\rangle-4\langle d\phi(\partial_r),\bar\nabla_{\partial_r}\tau(\phi)\rangle\big)d\sigma.
\end{align*}
In the final step we use 
\begin{align*}
(2-m)\int_{B_r}|d\phi|^2d\mu+r\int_{\partial B_r}|d\phi|^2d\sigma&=r^{m-1}\frac{d}{dr}r^{2-m}\int_{B_r}|d\phi|^2d\mu, \\
(4-m)\int_{B_r}|\tau(\phi)|^2d\mu+r\int_{\partial B_r}|d\phi|^2d\sigma&=r^{m-3}\frac{d}{dr}r^{4-m}\int_{B_r}|\tau(\phi)|^2d\mu,
\end{align*}
which is a direct consequence of the coarea-formula.
\end{proof}

It is obvious that the second to last term on the right hand side of \eqref{identity-general-rm}
is an obstacle when trying to derive a monotonicity formula.
Fortunately, for a certain class of interpolating sesqui-harmonic maps this contribution
vanishes and we can give a kind of monotonicity formula.

\begin{Satz}
Let \(\phi\colon\R^m\to N\) be a smooth solution of \eqref{euler-lagrange}
for which \(\tau(\phi)\) is orthogonal to the image of the map.
In addition, assume that \(m>4\) and by \(B_r\) we denote the geodesic ball of radius \(r\) in \(\R^m\).
Then the following monotonicity formula holds
\begin{align}
\label{monotonicity-formula}
R_1^{4-m}\int_{B_{R_1}}(\delta_1|d\phi|^2+&\delta_2|\tau(\phi)|^2+4\sqrt{m}|\delta_2||\bar\nabla d\phi|^2)d\mu \\
\nonumber&\leq R_2^{4-m}\int_{B_{R_2}}(\delta_1|d\phi|^2+\delta_2|\tau(\phi)|^2+4\sqrt{m}|\delta_2||\bar\nabla d\phi|^2)d\mu,
\end{align}
where \(0<R_1<R_2\).
\end{Satz}
\begin{proof}
By assumption, we have 
\begin{align*}
\langle\tau(\phi),d\phi(X)\rangle=0 
\end{align*}
for all vector fields \(X\) on \(M\).
Hence, equation \eqref{identity-general-rm} yields the following inequality
\begin{align*}
r^{m-2}\frac{d}{dr}r^{2-m}\delta_1\int_{B_r}|d\phi|^2d\mu+r^{m-4}\frac{d}{dr}r^{4-m}\delta_2\int_{B_r}|\tau(\phi)|^2d\mu
\geq-4\sqrt{m}|\delta_2|\int_{\partial B_r}|\bar\nabla d\phi|^2 d\sigma.
\end{align*}
Multiplying by \(r^{4-m}\) and integrating with respect to \(r\) from \(R_1\) to \(R_2\) we find
\begin{align*}
R_1^{4-m}\delta_2\int_{B_{R_1}}|\tau(\phi)|^2d\mu\leq &R_2^{4-m}\delta_2\int_{B_{R_2}}|\tau(\phi)|^2d\mu
+\int_{R_1}^{R_2}\big(r^{2}\frac{d}{dr}r^{2-m}\delta_1\int_{B_r}|d\phi|^2d\mu\big)dr\\
&+4\sqrt{m}|\delta_2|\int_{R_1}^{R_2}\big(r^{4-m}\int_{\partial B_r}|\bar\nabla d\phi|^2 d\sigma\big)dr.
\end{align*}
Using integration by parts we find
\begin{align*}
\int_{R_1}^{R_2}\big(r^{2}\frac{d}{dr}r^{2-m}\int_{B_r}|d\phi|^2d\mu\big)dr=&-2\int_{R_1}^{R_2}\big(r^{3-m}\int_{B_r}|d\phi|^2d\mu\big)dr \\
&+R_2^{4-m}\int_{B_{R_2}}|d\phi|^2d\mu-R_1^{4-m}\int_{B_{R_1}}|d\phi|^2d\mu.
\end{align*}
In addition, we find by the same arguments 
\begin{align*}
\int_{R_1}^{R_2}\big(r^{4-m}\int_{\partial B_r}|\bar\nabla d\phi|^2 d\sigma\big)dr=&
-(4-m)\int_{R_1}^{R_2}\big(r^{3-m}\int_{B_r}|\bar\nabla d\phi|^2 d\mu\big)dr \\
&+R_2^{4-m}\int_{B_{R_2}}|\bar\nabla d\phi|^2d\mu
-R_1^{4-m}\int_{B_{R_1}}|\bar\nabla d\phi|^2d\mu
\end{align*}
and the claim follows due to the assumption \(m>4\).
\end{proof}

\begin{Bem}
\begin{enumerate}
\item It is straightforward to generalize \eqref{monotonicity-formula} to the case that the domain is a Riemannian manifold.
\item The assumptions of the Theorem hold in the case that \(\phi\) is an isometric immersion.
However, in this case we also have \(|d\phi|^2=m\) such that \eqref{monotonicity-formula}
does not contain much information.
\end{enumerate}
\end{Bem}

\par\medskip
\textbf{Acknowledgements:}
The author gratefully acknowledges the support of the Austrian Science Fund (FWF) 
through the project P 30749--N35 ``Geometric variational problems from string theory''.

\bibliographystyle{plain}
\bibliography{mybib}
\end{document}